\documentclass[pdftex,10pt,a4paper,oneside]{article}

\usepackage[english]{babel}
\usepackage[utf8]{inputenc}
\usepackage[T1]{fontenc}

\usepackage{lmodern}

\usepackage{amsmath}
\usepackage{amsthm}

\usepackage{geometry}

\usepackage{authblk}

\usepackage[runin]{abstract}

\usepackage{titlesec}

\titleformat{\section}[block]{\large\bfseries}{\thesection}{1.0em}{}
\titleformat{\subsection}[runin]{\normalsize\bfseries}{\thesubsection}{0.5em}{}
\titleformat{\paragraph}[runin]{\normalsize\itshape}{\theparagraph}{}{}

\usepackage[section]{tocbibind}

\numberwithin{equation}{subsection}

\newcommand{\PAR}[1]{\subsection{#1\texorpdfstring{.}{}}}

\newenvironment{THEOREM}[1]{\subsection{Theorem%
\if #1\@\texorpdfstring{.}{} \else \ ({#1})\texorpdfstring{.}{} \fi}\it}{}

\newenvironment{LEMMA}[1]{\subsection{Lemma%
\if #1\@\texorpdfstring{.}{} \else \ ({#1})\texorpdfstring{.}{} \fi}\it}{}

\newcommand{\XELerr}[2]{(\ref{eq:#1})}      
\newcommand{\XEL}[2]{(\ref{eq:#1\if #2\@\else:#2\fi})}
\newcommand{\XEGerr}[2]{\ref{par:#1}(\ref{eq:#1})} 
\newcommand{\XEG}[2]{\ref{par:#1}(\ref{eq:#1\if #2\@\else:#2\fi})}
\newcommand{\XS}[1]{Sect.~\ref{sc:#1}}      
\newcommand{\XP}[1]{Par.~\ref{par:#1}}      
\newcommand{\XT}[1]{Thm.~\ref{cl:#1}}       
\newcommand{\XL}[1]{Lemma~\ref{cl:#1}}      

\def\PP{\hbox{\rlap{\boldmath$,$}\lower1.2pt\hbox{\boldmath$\cdot$}\,}}
\newcommand{\SM}{{\cal M}}
\newcommand{\mRW}{\Rightarrow}
\newcommand{\mLW}{\Leftarrow}

\newcommand{\RW}{\rightarrow}

\newcommand{\LRW}{\leftrightarrow}

\newcommand{\Za}{\BO{0}}
\newcommand{\Zb}{\tilde{\BO{0}}}
\newcommand{\Zc}{\hat{\BO{0}}}

\newcommand{\SN}{{\cal N}}
\newcommand{\MN}{{\cal M}_{\alpha}}

\newcommand{\Apc}{{\cal A}}

\newcommand{\BO}[1]{\mbox{\boldmath $#1$}}
\newcommand{\IT}[1]{{\it #1\/}}

\usepackage[pdftex,bookmarksnumbered,pagebackref]{hyperref} 

\hypersetup{
 pdftitle={On Herbrand Skeletons},
 pdfauthor={Paul J. Voda and Ján Komara},
 pdfsubject={Mathematical Logic},
 pdfkeywords={
  Herbrand's theorem,
  Herbrand skeletons
 },
}

\hypersetup{
 unicode,
 colorlinks=true,
 allcolors=[rgb]{0,0.5,0.5},
 urlcolor=[rgb]{0.58,0.0,0.83},
 breaklinks=true,
}

\begin{document}


\title{On Herbrand Skeletons}

\author{Paul J. Voda}
\author{Ján Komara}

\affil{
Institute of Informatics, Faculty of Mathematics and Physics,
\\
Comenius University, Mlynská dolina, 842 15 Bratislava, Slovakia.
\\
Technical Report mff-ii-02-1995
\\
July 20, 1995; revised January 26, 1996
\\
\{voda,komara\}@fmph.uniba.sk
}

\date{}

\maketitle


\renewcommand{\abstractname}{Abstract.}

\begin{abstract}
Herbrand's theorem plays an important role both in proof theory and in
computer science. Given a Herbrand skeleton, which is basically a
number specifying the count of disjunctions of the matrix, we would
like to get a computable bound on the size of terms which make the
disjunction into a quasitautology. This is an important problem in
logic, specifically in the complexity of proofs.  In computer
science, specifically in automated theorem proving, one hopes for an
algorithm which avoids the guesses of existential substitution axioms
involved in proving a theorem. Herbrand's theorem forms the very
basis of automated theorem proving where for a given number $n$ we
would like to have an algorithm which finds the terms in the $n$
disjunctions of matrices solely from the shape of the matrix. The
main result of this paper is that both problems have negative
solutions.
\end{abstract}


\section{Introduction}\label{sc:intro}

By the theorem of Herbrand we have for a quantifier-free $\phi$:
\[
\models \exists \bar{\BO{x}}\,\phi(\bar{\BO{x}})\
\mbox{iff}\
\models \phi(\bar{\BO{a}}_1) \vee \phi(\bar{\BO{a}}_2) \vee \cdots \vee
        \phi(\bar{\BO{a}}_n)
\]
for certain $n$ and sequences of terms
\( \bar{\BO{a}}_1, \ldots, \bar{\BO{a}}_n \).
The question whether for a given $n$ we can find such terms $\bar{\BO{a}}_i$
for arbitrary formulas $\phi$ is the problem of {\em Herbrand skeletons of size
$n$\/}. The terms $\bar{\BO{a}}_i$ {\em solve\/} the skeleton.

The problem is very important both in logic where we enquire about bounds on
the size of such terms
\cite{qq:hajekpudlak}
and also in automated theorem proving
(ATP) where one asks whether there is an algorithm which finds the terms from
$\phi$ without a guess involved in the choice of existential substitution
axioms.  This use of Herbrand's theorem is actually the very foundation of ATP.
This is because the theorem of Herbrand has influenced the oldest ATP proof
procedure of {\em resolution\/} \cite {qq:robinson} and Herbrand skeletons
directly appear in the modern ATP procedures based on the {\em connection
method\/} (in the form of the {\em multiplicity\/} of formulas)
\cite{qq:bibel:82} and on semantic tableaux (in the form of {\em free
variables\/}) \cite{qq:fitting}.

It was already known to Herbrand \cite{qq:herbrand} that when the
formula $\phi$ does not contain the identity $\doteq$ then we can
both effectively find the bounds and that there is an algorithm for
finding the terms (unification).
With identities permitted in the formulas $\phi$,
both problems were open although a recent result by Degtyarev and
Voronkov on undecidability of the so called {\em simultaneous rigid
E-unification\/} \cite{qq:degvor:rigid:1} can be used to show the
undecidability for the case when $n = 1$ (see
\XP{sk:rigid} for more details). However, this result does not readily extend
to the case $n > 1$.

Inspired by \cite{qq:degvor:rigid:1} and \cite{qq:krajpudlak} we settle in this paper
both problems negatively: for no $n$ there is a computable function in $\phi$
giving a bound on the size of a solution, nor there is an algorithm with input
$\phi$ for finding a solution.  We do it by the reduction of the celebrated
result of Matiyasevich on unsolvability of Diophantine equations \cite{qq:mat}
to the solvability of $n$-skeletons. In order to emphasize the purely logical
character of the problem of Herbrand skeletons we deliberately refrain from
using any specialized ATP and/or term-rewriting terminology.

In \XS{pc} we give our notation, \XS{sk} introduces the problem of
Herbrand skeletons. Section \ref{sc:sk1} deals with the undecidability
of Herbrand skeletons for $n = 1$, \XS{arpc} proves technical lemmas
needed for this.
Section \ref{sc:skn} presents the main result for $n \geq 1$, \XS{main}
proves technical lemmas for this.


\section{Notation and Logical Background}\label{sc:pc}

\PAR{Language of predicate calculus}\label{par:pc:lang}
The {\em language of first-order predicate calculus with identity} consists of
denumerably many variables, function and predicate symbols of all arities
(function symbols of arity $0$ are constants, predicate symbols of arity $0$
are propositional constants). We use $\BO{f}$, and $\BO{p}$ as metavariables
ranging over function and predicate symbols respectively,
\( \BO{x}, \BO{y}, \ldots \),
possibly with subscripts, as metavariables ranging over variables.

{\em Semiterms\/} are either variables or expressions
\( \BO{f}(\BO{a}_1,\BO{a}_2,\ldots,\BO{a}_n) \)
where $\BO{f}$ is a function symbol of arity $n$ and $\BO{a}_i$ are previously
constructed semiterms.  We use
\( \BO{a}, \BO{b}, \ldots \)
as metavariables ranging over semiterms, $\BO{k}$ as metavariables ranging
over constant symbols, and we write {\em constant\/} semiterms $\BO{k}()$ as
$\BO{k}$.

Any metavariable, say \BO{a}, written as $\bar{\BO{a}},$ denotes a possibly
empty sequence of objects (terms, constants, variables) denoted by the
metavariable. Sequences of variables $\bar{\BO{x}}$, $\bar{\BO{y}}$ used in a
certain context are always assumed to consist of pairwise distinct variables
without having any variables in common.

{\em Semiformulas\/} are constructed from {\em atomic\/} semiformulas
of identity
\( \BO{a} \doteq \BO{b} \)
and predicate applications
\( \BO{p}(\BO{a}_1,\BO{a}_2,\ldots,\BO{a}_n) \)
by {\em propositional\/} connectives $\neg$, $\vee$, $\wedge$, $\RW$ and {\em
quantifiers\/} $\exists$, $\forall$ in the usual way.  We use lowercase Greek
letters
\( \phi, \psi, \ldots \)
to range over semiformulas. Free and bound variables of semiformulas are
defined as usual.

{\em Terms\/} are semiterms without variables and {\em formulas\/} are
semiformulas without free variables.  As usual, a {\em structure} $\SM$ is
given by a non-empty domain and an interpretation $I$ of function and predicate
symbols.  We write $\BO{a}^I$ for the denotation of the term $\BO{a}$.  For a
formula $\phi$ we write
\( \SM \models \phi \)
to assert that $\phi$ is true in the structure $\SM$, $\models \phi$
means that $\phi$ is valid, i.e. that it holds in all structures.

We assume that the objects of first-order predicate calculus are
encoded into natural numbers in some of the usual ways.  Thus, for instance,
the set of terms is a subset of natural numbers.  By having variables,
semiterms, and semiformulas as natural numbers the metatheoretic predicates, as
for instance, $P(\phi)$ are numeric predicates. We may then say that $P(\phi)$
is a recursively enumerable (r.e.) predicate without having to speak of codes
of objects. We are also able to write $\BO{a} = \BO{b}$ to assert that both
terms are identical (note that atomic identity formulas, which are numbers, are
written as $\BO{a} \doteq \BO{b}$).

\PAR{Quasitautologies}\label{par:pc:qtaut}
{\em Quasitautologies\/} are quantifier-free formulas which are
tau\-to\-lo\-gi\-cal (propositional) consequences of {\em identity axioms\/}:
\[
\BO{k} \doteq \BO{k}
\]
\[
\BO{a} \doteq \BO{b} \RW \BO{b} \doteq \BO{a}
\]
\[
\BO{a} \doteq \BO{b} \wedge \BO{b} \doteq \BO{c} \RW \BO{a} \doteq \BO{c}
\]
\[
\BO{a}_1 \doteq \BO{b_1} \wedge \cdots \wedge \BO{a}_n \doteq \BO{b}_n \RW
        \BO{f}(\BO{a}_1,\ldots,\BO{a}_n) \doteq \BO{f}(\BO{b}_1,\ldots,\BO{b}_n)
\]
\[
\BO{a}_1 \doteq \BO{b_1} \wedge \cdots \wedge \BO{a}_n \doteq \BO{b}_n \wedge
        \BO{p}(\BO{a}_1,\ldots,\BO{a}_n) \RW
        \BO{p}(\BO{b}_1,\ldots,\BO{b}_n) \, .
\]
For a quantifier-free formula $\phi$ we have
\[
\mbox{$\models \phi$ iff $\phi$ is a quasitautology.}
\]
The predicate of being a quasitautology is
primitive recursive (see, for instance, \cite{qq:boolos}).

\PAR{Language of arithmetic}\label{par:pc:arithm}
The first-order {\em language of \em arithmetic\/} $L$ consists of the constant
$\Za$, unary function symbol $S$, and two binary function symbols $+$, and
`$\cdot$'. The structure
\[
\SN = \langle N, \Za, S, +, \cdot \rangle
\]
is over the domain of natural numbers with the standard interpretation of the
constant $\Za$ as zero, $S$ as the successor function, $+$ as addition, and
`$\cdot$' as multiplication.

For every constant \BO{k} we define {\em $\BO{k}$-numerals\/} as terms of the
form $S^m(\BO{k})$ for some $m \geq 0$. Here
\( S^0(\BO{a}) = \BO{a} \)
and
\( S^{m+1}(\BO{a}) = S\,S^m(\BO{a}) \).

The class of {\em diophantine semiformulas\/} is composed from atomic
semiformulas
\( \BO{a} + \BO{b} \doteq \BO{c} \)
and
\( \BO{a} \cdot \BO{b} \doteq \BO{c} \)
by conjunctions. Here the terms $\BO{a}$, $\BO{b}$, and $\BO{c}$ are either
variables or $\Za$-numerals.

By the theorem of Matiyasevich~\cite{qq:mat} every recursively enumerable
predicate $R(m)$ can be represented by a diophantine semiformula
$\psi(\BO{x},\bar{\BO{x}})$ with only the indicated variables
free such that for every number $m$
\begin{equation}\label{eq:pc:arithm}
R(m)\
\mbox{iff}\
\SN \models \exists \bar{\BO{x}}\, \psi(S^m(\Za),\bar{\BO{x}})\, .
\end{equation}
Although the predicate $\SN \models \psi$ restricted to diophantine formulas is
primitive recursive, the same predicate restricted to existentially closed
diophantine semiformulas is only recursively enumerable.


\section{Herbrand Skeletons}\label{sc:sk}

Herbrand's theorem is the formal basis for ATP.  For
our purposes it is convenient to state it in terms of solvability of formulas.

\PAR{Solvability of formulas}\label{par:sk:solv}
Among the unlimited supply of constant symbols in the language of predicate
calculus we single out the constants
\( \ast, \ast_1, \ast_2, \ldots \)
and call them {\em unknowns\/}.  We use the metavariables
\( \BO{\ast}, \BO{\ast}_1, \ldots \)
to range over unknowns. Sequences of unknowns are denoted by $\bar{\BO{\ast}}$
with the same conventions as for sequences of variables.

We will indicate by writing $\phi(\bar{\BO{\ast}})$ that the unknowns occurring
in the formula $\phi$ are among $\bar{\BO{\ast}}$. For a sequence of terms
$\bar{\BO{a}}$ of the same length we will write $\phi(\bar{\BO{a}})$ for the
formula obtained by the simultaneous replacement in $\phi$ of unknowns by the
respective terms in $\bar{\BO{a}}$.

Terms $\bar{\BO{a}}$ (without unknowns) are called a {\em solution\/} of the
formula $\phi(\bar{\BO{\ast}})$ if $\models \phi(\bar{\BO{a}})$. A formula
$\phi$ is {\em solvable\/} if it has a solution.

We mention some obvious facts about solvability.  If $\phi \wedge \psi$ is
solvable then both $\phi$ and $\psi$ are.
If $\phi$ and $\psi$ are solvable and they do not share unknowns then
also $\phi \wedge \psi$ is solvable.
If $\phi$ or $\psi$ is solvable
then also $\phi \vee \psi$ is. The converse does not hold even if
the unknowns are not shared.

\PAR{Herbrand's theorem}\label{par:sk:herb}
A formula of the form
\( \exists \bar{\BO{x}}\, \phi(\bar{\BO{x}}) \)
with $\phi$ quantifier-free and without unknowns is called {\em existential
formula\/}. The semiformula $\phi$ is its {\em matrice\/}.

The first part of the theorem of Herbrand says that for every formula $\psi$ we
can find an existential formula
\( \exists \bar{\BO{x}}\,\phi(\bar{\BO{x}}) \)
such that
\[
\models \psi\ \mbox{iff}\ \models \exists \bar{\BO{x}}\, \phi(\bar{\BO{x}})\, .
\]
This part can proved by the elimination of universal quantifiers from $\psi$ by
means of Skolem functions (see, for instance, \cite{qq:shoenfield}).

For an existential formula $\psi$ with the matrix $\phi(\bar{\BO{x}})$
we call any formula
\begin{equation}\label{eq:sk:herb:1}
\phi(\bar{\BO{\ast}}_1) \vee \phi(\bar{\BO{\ast}}_2) \vee \cdots \vee
\phi(\bar{\BO{\ast}}_n)
\end{equation}
a {\em Herbrand skeleton of $\psi$ of size $n$\/}.

The second part of Herbrand's theorem says that for an
existential formula $\psi$
we have
\[
\models \psi\ \mbox{iff}\
\psi\ \mbox{has a solvable skeleton of some size $n$.}
\]

We note that skeletons of the same size differ only in the names of unknowns.
Hence, the solvability of any one of them implies the solvability of all
skeletons of the same (and larger) size.

Herbrand skeletons are quantifier-free formulas. Thus the test whether a given
sequence of terms is a solution to a given skeleton involves the primitive
recursive test whether the formula obtained from the skeleton by the
replacement of unknowns by the terms is a quasitautology.

We define the predicate $\IT{Sk}(n,\psi)$ as
\begin{quote}
$\IT{Sk}(n,\psi)$ iff $\psi$ is an existential formula with
a solvable skeleton of size $n$.
\end{quote}
Hence, for an existential formula $\psi$ we have
$\models \psi$ iff
$\IT{Sk}(n,\psi)$ for some $n$.
For every number $n$ we define the predicate $\IT{Sk}_n(\psi)$ as
\begin{quote}
$\IT{Sk}_n(\psi)$ iff $\IT{Sk}(n,\psi)$.
\end{quote}

\PAR{Solvability of Herbrand skeletons}\label{par:sk:skel}
Predicate calculus is semi-de\-ci\-da\-ble, i.e.  the
predicate $\models \psi$ is recursively enumerable but not recursive. By the
first part of Herbrand's theorem neither the restriction of $\models \psi$ to
existential formulas is recursive. Hence, the equivalent predicate in $\psi$:
\begin{quote}
$\IT{Sk}(n,\psi)$ for some $n$
\end{quote}
is not recursive although it is recursively enumerable.

The recursive enumerability of this predicate has two `degrees of freedom' as
it involves two guesses: first the number $n$ and then the solution.  It was
hoped in ATP circles that the second guess was not needed and that there was an
algorithm which would find a solution or would determine that a formula is
unsolvable.  In other words, it was hoped that the predicates $\IT{Sk}_n(\psi)$
were recursive. The main result of this paper is that this is not the case for
any $n \geq 1$.

\PAR{Example}\label{par:sk:ex1}
Consider an existential formula
\(
\exists x\, (p(a) \vee p(b) \RW p(x))
\)
where $a$ and $b$ are different constants. We can solve its $2$-skeleton
\[ (p(a) \vee p(b) \RW p(\ast)) \vee (p(a) \vee p(b) \RW p(\ast_1)) \]
since
\begin{equation}\label{eq:sk:ex1:2}
\models (p(a) \vee p(b) \RW p(c)) \vee (p(a) \vee p(b) \RW p(c))\, .
\end{equation}
On the other hand, the $1$-skeletons are not solvable.

\PAR{Simultaneous rigid E-unification}\label{par:sk:rigid}
{\em A simultaneous rigid $E$-unification problem\/} ({\em SREU problem\/} for
short) is a problem of finding a solution of a formula $\phi(\bar{\BO{\ast}})$
which is a conjunction of $n$-formulas of the form
\begin{equation}\label{eq:sk:rigid:1}
\BO{a}_1 \doteq \BO{b}_1 \wedge \cdots \wedge \BO{a}_m \doteq \BO{b}_m \RW
\BO{a} \doteq \BO{b}\, .
\end{equation}
We do not exclude $m = 0$ in which case the above formula is just the identity
$\BO{a} \doteq \BO{b}$. For $n=1$ we have a {\em rigid E-unification
problem\/}.

It has been known for long that (non-simultaneous) rigid E-unification is
decidable (see \cite{qq:GNPS:88},\cite{qq:GNPS:90}; for a more elementary proof
see \cite{qq:kogel}) while the decidability of SREU
has been an open problem. It was
recently settled negatively by Degtyarev and Voronkov \cite{qq:degvor:rigid:1}.

The reader will note that the problem of finding a solution of
$\phi(\bar{\BO{\ast}})$ is equivalent to the problem of whether the
existential formula
\( \exists \bar{\BO{x}}\, \phi(\bar{\BO{x}}) \)
has a solvable $1$-skeleton. Thus the undecidability of SREU
implies the undecidability of $\IT{Sk}_1(\psi)$.  However, the
undecidability of $\IT{Sk}_n(\psi)$ for any $n$, which is our main result given
in \XT{skn:col}, is not a direct consequence.

We will outline in \XP{sk:r} a procedure,
which is a matter of folklore in ATP circles.
The procedure converts a Herbrand skeleton of
size $n$ to a finite class of SREU
problems such that the skeleton is solvable iff at least one
SREU from the class is.

When all formulas of a SREU problem
are identities, the problem becomes a {\em (syntactical) unification
problem\/}. That this is decidable was already known  to Herbrand
\cite{qq:herbrand}.
Consequently, when an existential formula does not contain the
identity $\doteq$ then for any skeleton of size $n$ the conversion
procedure yields a finite number of unification problems (i.e.
$m=0$ in each of the problems). Hence, $\IT{Sk}_n(\psi)$ restricted to such
existential formulas is decidable. This has been also known to Herbrand, see
also
\cite{qq:buss:intro,qq:buss:herb}.

\PAR{Converting Herbrand skeletons to SREU problems}\label{par:sk:r}
In this paragraph we reduce in the following sense the problem of solvability
of Herbrand skeletons to a class of SREU problems:
\begin{itemize}
\item[]
To every quantifier-free formula
\( \phi(\bar{\BO{\ast}}) \)
we can primitively
recursively find a finite class $\Gamma$ of SREU problems
which are {\em solution equivalent\/} in the sense that
every solution to
\( \phi(\bar{\BO{\ast}}) \)
solves at least one problem from $\Gamma$
and vice versa, every solution of a problem from $\Gamma$ solves
\( \phi(\bar{\BO{\ast}}) \).
\end{itemize}
Let 
\( \phi(\bar{\BO{\ast}}) \)
be a quantifier-free formula. The transformation consists of three steps.

\paragraph*{(i): Transformation to conjunction of clauses.}
We first convert
\( \phi(\bar{\BO{\ast}}) \)
into an equivalent conjunction of clauses
$\phi_1(\bar{\BO{\ast}})$.
A {\em clause} is of a form
\begin{equation}\label{eq:sk:r:1}
A_1 \wedge \cdots \wedge A_n \RW B_1 \vee \cdots \vee B_m
\end{equation}
for atomic formulas
\( A_1, \ldots, A_n, B_1, \ldots, B_m \).
For $m=1$ the clause is called a {\em Horn\/} clause.  We do not exclude the
case when $n=0$ in which case \XEL{sk:r}{1} stands for
\(  B_1 \vee \cdots \vee B_m \)
or the case when $m=0$ in which case we put the formula $c \doteq d$ for two
{\em new\/} distinct constants in the consequent (body) of the clause.

We set $\Gamma_0 = \{ \phi_1(\bar{\BO{\ast}}) \}$ and observe
that $\Gamma_0$ and
\( \phi(\bar{\BO{\ast}}) \)
are solution equivalent.

\paragraph{(ii): Transformation to Horn clauses.}
Assume that we are given a finite class $\Gamma_i$ of formulas which are
conjunctions of clauses such that $\Gamma_i$
is solution equivalent to
\( \phi(\bar{\BO{\ast}}) \).
If some formula $\phi'(\bar{\BO{\ast}}) \in \Gamma_i$
has a form
\[
-\,-\,- \wedge
(A_1 \wedge \cdots \wedge A_n \RW B_1 \vee \cdots \vee B_m)
\wedge -\,-\,-
\]
with $m > 1$ we replace the formula in $\Gamma_i$ by the set of formulas
\[ -\,-\,- \wedge
   (A_1 \wedge \cdots \wedge A_n \RW B_j) \wedge -\,-\,-
\]
for $1 \leq j \leq m$. We obtain a new class $\Gamma_{i+1}$ which
is solution equivalent to
\( \phi(\bar{\BO{\ast}}) \).
We repeat the process as long as $\Gamma_{i+i_1}$ contains non Horn-clauses.

\paragraph{(iii): Elimination of predicate symbols.}
Assume that we are given a finite class $\Gamma_k$ of formulas which are
conjunctions of Horn clauses such that $\Gamma_k$
is solution equivalent to
\( \phi(\bar{\BO{\ast}}) \).
If some formula $\phi'(\bar{\BO{\ast}}) \in \Gamma_k$
is not a SREU problem then one of the following cases
must obtain:
\begin{itemize}
\item
The formula $\phi'(\bar{\BO{\ast}})$ has a form
\[ -\,-\,- \wedge
   (\cdots \RW \BO{p}(\bar{\BO{a}})) \wedge -\,-\,-
\]
and the predicate symbol $\BO{p}$ does not occur in the antecedent of the
clause.  Then the formula is unsolvable as it can be always falsified in a
suitable structure. We delete the formula from $\Gamma_k$.
\item
The formula $\phi'(\bar{\BO{\ast}})$ has a form
\[
-\,-\,- \wedge
(\cdots \wedge \BO{q}(\bar{\BO{b}}) \wedge \cdots \RW \BO{p}(\bar{\BO{a}}))
\wedge -\,-\,-
\]
where \BO{q} and \BO{p} are distinct
predicate symbols. Then we replace the formula in $\Gamma_k$ by the formula
\[
-\,-\,- \wedge
( \cdots \wedge \cdots \RW \BO{p}(\bar{\BO{a}}))
\wedge -\,-\,- \, .
\]
\item
The formula $\phi'(\bar{\BO{\ast}})$ has a form
\[
-\,-\,- \wedge
(\cdots \wedge \BO{p}(\BO{b}_1,\ldots, \BO{b}_n) \wedge \cdots \RW
	\BO{p}(\BO{a}_1,\ldots,\BO{a}_n))
\wedge -\,-\,-\, .
\]
Then we replace the formula in $\Gamma_k$ by two formulas
\[
-\,-\,- \wedge
(\cdots \wedge \cdots \RW \BO{p}(\BO{a}_1,\ldots,\BO{a}_n)))
\wedge -\,-\,- \, ,
\]
\[
-\,-\,- \wedge
(\cdots \wedge \cdots \RW \BO{a}_1 \doteq \BO{b}_1)
	\wedge \ldots \wedge
(\cdots \wedge \cdots \RW \BO{a}_n \doteq \BO{b}_n)
\wedge -\,-\,- \ \, .
\]
\end{itemize}
By the above changes we obtain a new class $\Gamma_{k+1}$ which is
solution equivalent to
\( \phi(\bar{\BO{\ast}}) \).
We repeat the process as long as $\Gamma_{k+k_1}$ contains
formulas which are not SREU's.

\PAR{Example}\label{par:sk:ex2}
The above conversion is demonstrated with the
formula
\[
\exists x\, (p(a) \wedge p(b) \wedge (x \doteq a \vee x \doteq b) \RW p(c)) \, .
\]
Here $a$, $b$, and $c$ are different constants.
Note that the constant $c$ is a
solution of its $1$-skeleton
$p(a) \wedge p(b) \wedge (\ast \doteq a \vee \ast \doteq b) \RW p(c)$.
Converting the skeleton to a conjunction of clauses yields
\[
(p(a) \wedge p(b) \wedge \ast \doteq a \RW p(c)) \wedge
(p(a) \wedge p(b) \wedge \ast \doteq b \RW p(c)) \, .
\]
Elimination of predicate symbols leads to four
SREU problems:
\begin{eqnarray}
(\ast \doteq a \RW a \doteq c) \wedge (\ast \doteq b \RW a \doteq c)\, ,
\nonumber
\\
(\ast \doteq a \RW a \doteq c) \wedge (\ast \doteq b \RW b \doteq c)\, ,
\label{eq:sk:ex2:42}
\\
(\ast \doteq a \RW b \doteq c) \wedge (\ast \doteq b \RW a \doteq c)\, ,
\nonumber
\\
(\ast \doteq a \RW b \doteq c) \wedge (\ast \doteq b \RW b \doteq c)\, .
\nonumber
\end{eqnarray}
Problem \XEL{sk:ex2}{42} is the only solvable one as it is solved by the
constant $c$.


\section{Non-recursiveness of \texorpdfstring{$Sk_1(\psi)$}{Sk₁(Ψ)}}\label{sc:sk1}

\PAR{Language of arithmetic formulas in predicate calculus}\label{sc:arpc:lang}
We wish to simulate arithmetic by certain quantifier-free semiformulas of
predicate calculus. The semiformulas will be in the language $P$ consisting of
constants $\Za$, $\Zc$, $\Zb$, $k$, $\tilde{k}$, of the unary function symbol
$S$, and of the binary function symbol `\PP'.  We write the binary symbol in
the infix form $\BO{a}\PP\BO{b}$ where $\BO{a}\PP\BO{b}\PP\BO{c}$ associates to
the right, i.e. it is read as $\BO{a}\PP(\BO{b}\PP\BO{c})$.  The function
symbol $S$ will simulate the successor function, while the function symbol
`\PP' will play the role of a pairing function (\IT{cons} of LISP).

We will define in Paragraphs \ref{par:arpc:num:df}, \ref{par:arpc:add:df}, and
\ref{par:arpc:mul:df} quantifier-free semiformulas
\( \IT{Num}(x) \),
\( \IT{Add}(x,y,z,w) \),
and
\( \IT{Mul}(x,y,z,w,\tilde{w}) \)
of the language $P$ with all of their free variables indicated.  The
semiformulas simulate arithmetic in predicate calculus as can be seen from the
following lemma which will be proved in \XS{arpc}.
\begin{LEMMA}{}\label{cl:arpc}
\begin{itemize}
\item[\rm (a)]
$\IT{Num}(\BO{\ast})$ is solved exactly by \Za-numerals,
\item[\rm (b)]
\( \IT{Add}(S^m(\Za),S^p(\Za),S^q(\Za),\BO{\ast}) \)
is solvable iff
$\SN \models S^m(\Za) + S^p(\Za) \doteq S^q(\Za)$,
\item[\rm (c)]
$\IT{Mul}(S^m(\Za),S^p(\Za),S^q(\Za),\BO{\ast},\BO{\ast}_1)$
is solvable iff
\( \SN \models S^m(\Za) \cdot S^p(\Za) \doteq S^q(\Za) \).
\end{itemize}
\end{LEMMA}

\PAR{PC-arithmetic semiformulas}\label{par:sk1:A}
We simulate arithmetic by certain quan\-ti\-fier-free semiformulas of the
language $P$ where we use two disjoint sets of variables:
\( x_1, x_2, x_3, \ldots \)
called {\em numeric\/} variables, and
\( w_1, w_2, w_3, \ldots \)
called {\em table\/} variables. We will use \BO{x} and \BO{w} as metavariables
ranging over numeric and table variables respectively.
The semiformulas are built up from the semiformulas
\( \IT{Num}(\BO{a}) \),
\( \IT{Add}(\BO{a},\BO{b},\BO{c},\BO{w}) \),
and
\( \IT{Mul}(\BO{a},\BO{b},\BO{c},\BO{w},\tilde{\BO{w}}) \)
by conjunctions. Here the terms \BO{a}, \BO{b}, and \BO{c} are either
$\Za$-numerals or numeric variables.

We associate with every diophantine semiformula $\psi(\bar{\BO{x}})$ a
quantifier-free semiformula $\phi(\bar{\BO{x}},\bar{\BO{w}})$ of $P$ called a
{\em PC-arithmetic semiformula\/}. The class of PC-arithmetic semiformulas is
denoted by $\Apc$.  The association is defined inductively as follows:
\begin{itemize}
\item
\( \BO{a} + \BO{b} \doteq \BO{c} \)
is associated with any semiformula of the form
\[
\IT{Num}(\BO{a}) \wedge \IT{Num}(\BO{b}) \wedge\IT{Num}(\BO{c}) \wedge
\IT{Add}(\BO{a},\BO{b},\BO{c},\BO{w})\, ,
\]

\item
\( \BO{a} \cdot \BO{b} \doteq \BO{c} \)
is associated with any semiformula of the form
\[
\IT{Num}(\BO{a}) \wedge \IT{Num}(\BO{b}) \wedge\IT{Num}(\BO{c}) \wedge
\IT{Mul}(\BO{a},\BO{b},\BO{c},\BO{w}_1,\BO{w}_2)\, ,
\]
\item
a diophantine semiformula
\( \psi_1(\bar{\BO{x}}) \wedge \psi_2(\bar{\BO{x}}) \)
is associated with any semiformula of the form
\[
\phi_1(\bar{\BO{x}},\bar{\BO{w}}_1) \wedge
\phi_2(\bar{\BO{x}},\bar{\BO{w}}_2)\, ,
\]
where
\( \psi_1(\bar{\BO{x}}) \)
and
\( \psi_2(\bar{\BO{x}}) \)
are associated with
\( \phi_1(\bar{\BO{x}},\bar{\BO{w}}_1) \)
and
\( \phi_2(\bar{\BO{x}},\bar{\BO{w}}_2) \)
respectively and the table variables $\bar{\BO{w}}_1$ and $\bar{\BO{w}}_2$
are disjoint.
\end{itemize}

It is easy to see that if the diophantine semiformula
\( \psi(\bar{\BO{x}}) \)
is associated with
\( \phi(\bar{\BO{x}},\bar{\BO{w}}) \in \Apc \)
then both semiformulas contain the same numeric variables
and every numeric variable \BO{x} of
\( \phi(\bar{\BO{x}},\bar{\BO{w}}) \)
occurs in a semiformula
\( \IT{Num}(\BO{x}) \).

\PAR{Invariancy of association under substitution}\label{par:sk1:sub}
We will use the following fact:
\begin{quote}
if the diophantine semiformula
\( \psi(\BO{x},\bar{\BO{x}}) \)
is associated with
\( \phi(\BO{x},\bar{\BO{x}},\bar{\BO{w}}) \in \Apc \)
then
\( \psi(S^m(\Za),\bar{\BO{x}}) \)
is associated with
\( \phi(S^m(\Za),\bar{\BO{x}},\bar{\BO{w}}) \in \Apc \)
\end{quote}
which is easily proved by induction on
\( \psi(\BO{x},\bar{\BO{x}}) \).

For the proof of the undecidability of 1-skeletons in \XT{sk1:col} we need some
auxiliary propositions.

\begin{LEMMA}{}\label{cl:sk1:cl1}\label{par:sk1:cl1}
If the diophantine formula $\psi$ is associated with
\( \phi(\bar{\BO{w}}) \in \Apc \) then
\( \SN \models \psi \)
iff
\( \phi(\bar{\BO{\ast}}) \) is solvable.
\end{LEMMA}

\begin{proof}
By induction on $\psi$. If $\psi$ is
\( S^m(\Za) + S^p(\Za) \doteq S^q(\Za) \)
then
\( \phi(\bar{\BO{\ast}})\)
has a form
\[
\IT{Num}(S^m(\Za)) \wedge \IT{Num}(S^p(\Za)) \wedge\IT{Num}(S^q(\Za)) \wedge
\IT{Add}(S^m(\Za),S^p(\Za),S^q(\Za),\BO{\ast})\; .
\]
The first three conjuncts are valid by \XL{arpc}(a) and the equivalence follows
directly from \XL{arpc}(b). The case when $\psi$ is
\( S^m(\Za) \cdot S^p(\Za) \doteq S^q(\Za) \)
is similar and uses \XL{arpc}(c). If $\psi$ is
\( \psi_1 \wedge \psi_2 \)
then
\( \phi(\bar{\BO{\ast}})\)
has a form
\( \phi_1(\bar{\BO{\ast}}_1) \wedge \phi_2(\bar{\BO{\ast}}_2) \)
where $\bar{\BO{\ast}}$ is partitioned into two disjoint sequences
\( \bar{\BO{\ast}}_1$ and $\bar{\BO{\ast}}_2 \).
Hence,
\( \SN \models \psi \)
iff
\( \SN \models \psi_1 \)
and
\( \SN \models \psi_2 \)
iff, by inductive hypotheses,
\( \phi_1(\bar{\BO{\ast}}_1) \)
and
\( \phi_2(\bar{\BO{\ast}}_2) \)
are solvable iff, because of disjointness of unknowns,
\( \phi_1(\bar{\BO{\ast}}_1) \wedge \phi_2(\bar{\BO{\ast}}_2) \)
is solvable.
\end{proof}

\begin{LEMMA}{}\label{cl:sk1:cl2}\label{par:sk1:cl2}
Let
\( \phi(\bar{\BO{x}},\bar{\BO{w}}) \in \Apc \)
be a semiformula with all of its numeric variables indicated. If the formula
\( \phi(\bar{\BO{a}},\bar{\BO{\ast}}) \)
is solvable then the terms $\bar{\BO{a}}$ are $\Za$-numerals.
\end{LEMMA}

\begin{proof}
We recall (\XP{sk1:A}) that every term \BO{a} of $\bar{\BO{a}}$
is substituted for $\BO{x}$ in some conjunct
\( \IT{Num}(\BO{x}) \)
of
\( \phi(\bar{\BO{x}},\bar{\BO{w}}) \).
Thus if \( \models \phi(\bar{\BO{a}},\bar{\BO{b}}) \)
for some terms $\bar{\BO{b}}$ then for every term \BO{a} of
$\bar{\BO{a}}$ we have
\( \models \IT{Num}(\BO{a}) \)
and by \XL{arpc}(a), the terms $\bar{\BO{a}}$ are $\Za$-numerals.
\end{proof}

\begin{LEMMA}{}\label{cl:sk1:cl3}\label{par:sk1:cl3}
Let the diophantine semiformula
\( \psi(\bar{\BO{x}}) \) with all of its variables indicated
be associated with some
\( \phi(\bar{\BO{x}},\bar{\BO{w}}) \in \Apc \).
Then
\( \SN \models \exists \bar{\BO{x}}\, \psi(\bar{\BO{x}}) \)
iff
\( \phi(\bar{\BO{\ast}},\bar{\BO{\ast}}_1) \)
is solvable.
\end{LEMMA}

\begin{proof}
The diophantine semiformula
\( \psi(\bar{\BO{x}}) \)
has the same set of numeric variables as
\( \phi(\bar{\BO{x}},\bar{\BO{w}}) \).
Thus
\( \SN \models \exists \bar{\BO{x}}\, \psi(\bar{\BO{x}}) \)
iff
\( \SN \models \psi(\bar{\BO{a}}) \)
for some $\Za$-numerals $\bar{\BO{a}}$ iff, by \XP{sk1:sub} and \XL{sk1:cl1},
\( \phi(\bar{\BO{a}},\bar{\BO{\ast}}_1) \)
is solvable for some $\Za$-numerals $\bar{\BO{a}}$ iff, by \XL{sk1:cl2} (in the
direction $\Leftarrow$),
\( \phi(\bar{\BO{\ast}},\bar{\BO{\ast}}_1) \)
is solvable.
\end{proof}

\begin{THEOREM}{}\label{cl:sk1:cl}\label{par:sk1:cl}
To every recursively enumerable predicate $R(m)$ there is a semiformula
\( \phi(\BO{x},\bar{\BO{x}},\bar{\BO{w}}) \in \Apc \)
such that for all $m$
\begin{equation}\label{eq:sk1:cl:1}
R(m)\ \mbox{iff}\
\IT{Sk}_1(\exists \bar{\BO{x}} \exists \bar{\BO{w}}\,
\phi(S^m(\Za),\bar{\BO{x}},\bar{\BO{w}}))\, .
\end{equation}
\end{THEOREM}

\begin{proof}
By the theorem of Matiyasevich there is a diophantine semiformula
\( \psi(\BO{x},\bar{\BO{x}}) \)
with all free variables indicated such that for every number $m$
\begin{equation}\label{eq:sk1:cl:2}
R(m)\ \mbox{iff}\
\SN \models \exists \bar{\BO{x}}\, \psi(S^m(\Za),\bar{\BO{x}})\, .
\end{equation}
Take a semiformula
\( \phi(\BO{x},\bar{\BO{x}},\bar{\BO{w}}) \in \Apc \)
associated to
\( \psi(\BO{x},\bar{\BO{x}}) \).
Then $R(m)$ holds iff, by \XEL{sk1:cl}{2},
\( \SN \models \exists \bar{\BO{x}}\, \psi(S^m(\Za),\bar{\BO{x}}) \)
iff, by \XP{sk1:sub} and \XL{sk1:cl3},
\( \phi(S^m(\Za),\bar{\BO{\ast}},\bar{\BO{\ast}}_1) \)
is solvable iff
\(
\IT{Sk}_1(\exists \bar{\BO{x}} \exists \bar{\BO{w}}\,
\phi(S^m(\Za),\bar{\BO{x}},\bar{\BO{w}}))
\).
\end{proof}

\begin{THEOREM}{}\label{cl:sk1:col}\label{par:sk1:col}
The predicate $\IT{Sk}_1(\psi)$ is not recursive.
\end{THEOREM}

\begin{proof}
Take any recursively enumerable but not recursive predicate $R(m)$ and obtain a
semiformula
\( \phi(\BO{x},\bar{\BO{x}},\bar{\BO{w}}) \in \Apc \)
from \XT{sk1:cl}. We can clearly find a primitive recursive function $f$ such
that
\[
f(m) = \exists \bar{\BO{x}} \exists \bar{\BO{w}}\,
	\phi(S^m(\Za),\bar{\BO{x}},\bar{\BO{w}})\, .
\]
Then $R(m)$ iff $\IT{Sk}_1(f(m))$ by \XEG{sk1:cl}{1}. If the predicate
$\IT{Sk}_1(\psi)$ were recursive so would be $R(m)$.
\end{proof}


\section{Simulation of Arithmetic}\label{sc:arpc}

In this section we will define the semiformulas \IT{Num}, \IT{Add}, and
\IT{Mul} simulating arithmetic in predicate calculus and prove \XL{arpc}. This
will finish the proof of the undecidability of $\IT{Sk}_1$ (\XT{sk1:col}).
The section is rather technical in that all proofs are carried out in detail.
We do this on purpose in order to demonstrate that the problem of Herbrand
skeletons is a purely logical problem albeit with extremely important
consequences for ATP.  Hence, we feel that the solution should be expressed in
the well-developed apparatus of predicate calculus (see for instance
\cite{qq:shoenfield}) without any detours through the terminology and
techniques of ATP and/or term rewriting.  We start with a lemma which is used
in ATP and term rewriting more or less automatically although its proof
requires non-trivial properties of predicate calculus.

\begin{LEMMA}{}\label{cl:pc:cl}
\begin{enumerate}
\item[\rm (a)]
For a semiformula $\phi(\BO{x})$ with at most \BO{x} free, term \BO{a}, and
constant $\BO{k}$ occurring neither in $\phi(\BO{x})$ nor in \BO{a} we have
\( \models \BO{k} \doteq \BO{a} \RW \phi(\BO{k}) \)
iff
\( \models \phi(\BO{a}) \).
\item[\rm (b)]
\( \models \BO{a} \doteq \BO{b} \)
iff
\( \BO{a} = \BO{b} \).
\end{enumerate}
\end{LEMMA}

\begin{proof}
(a):
We have
\( \models \BO{k} \doteq \BO{a} \RW \phi(\BO{k}) \)
iff, by the theorem on constants in \cite{qq:shoenfield},
\( \models \forall \BO{x} (\BO{x} \doteq \BO{a} \RW \phi(\BO{x})) \)
iff, by the third corollary of the equality theorem in \cite{qq:shoenfield},
\( \models \phi(\BO{a}) \).

(b): Clearly, if
\( \BO{a} = \BO{b} \)
then
\( \models \BO{a} \doteq \BO{b} \).
For the reverse direction consider a structure $\SM$ with the domain consisting
of all terms (which are a subset of natural numbers) and with the
interpretation $I$ of function symbols such that
\( \BO{f}^I(\bar d) = f(\bar d) \).
Clearly, $\BO{c}^I = \BO{c}$ for all terms.  If $\BO{a} \neq \BO{b}$ then
$\BO{a}^I \neq \BO{b}^I$, i.e.
\( \SM \models \BO{a} \not\doteq \BO{b} \),
and so
\( \not\models \BO{a} \doteq \BO{b} \).
\end{proof}

\PAR{Numerals}\label{par:arpc:num:df}
Denote by $\IT{Num}(x)$ the semiformula
\( \Za \doteq S(\Za) \RW \Za \doteq x \)
and by $\widetilde{\IT{Num}}(x)$ the semiformula
\( \Zb \doteq S(\Zb) \RW \Zb \doteq x \).

\begin{LEMMA}{}\label{cl:arpc:a}
\begin{itemize}
\item[\rm (a)]
$\IT{Num}(\BO{\ast})$ is solved exactly by \Za-numerals,
\item[\rm (b)]
$\widetilde{\IT{Num}}(\BO{\ast})$ is solved exactly by $\Zb$-numerals.
\end{itemize}
\end{LEMMA}

\begin{proof}
We prove only the part (a) as the proof of (b) is similar.
\( \models \IT{Num}(S^n(\Za)) \)
is proved by a straightforward induction on $n$. Conversely, if the term \BO{a}
is not a $\Za$-numeral then it must be the case that
\( \BO{a} = S^n(\BO{f}(\bar{\BO{b}})) \)
for a function symbol $\BO{f}$ different from $\Za$ and $S$, some number $n$,
and terms $\bar{\BO{b}}$. Consider a structure $\SM$ with the domain $\{0,1\}$
and the interpretation of function symbols $I$ such that $\Za^I = 0$, $S^I(d) =
d$ for all $d$ in the domain, and $\BO{g}^I(\bar d) = 1$, for all other symbols
and all $\bar d$ in the domain.  We clearly have
\( \Za^I = 0 = (S(\Za))^I \),
i.e.
\( \SM \models \Za \doteq S(\Za) \)
and also
\[
(S^n(\BO{f}(\bar{\BO{b}})))^I = (S^I)^n((\BO{f}(\bar{\BO{b}}))^I) =
(S^I)^n(1) = 1 \neq 0 = \Za^I
\, .
\]
Hence,
\(
\SM \not \models \Za \doteq S(\Za) \RW \Za \doteq S^n(\BO{f}(\bar{\BO{b}}))
\).
Thus
\( \not \models \IT{Num}(\BO{a}) \).
\end{proof}

\PAR{Proof of \XL{arpc}(a)}
This is \XL{arpc:a}(a). \qed

\PAR{Similar numerals}\label{par:arpc:sim:df}
Denote by $\IT{Sim}(x,y)$ the semiformula
$\Za \doteq \Zb \RW x \doteq y$.

\begin{LEMMA}{}\label{cl:arpc:sim}\label{par:arpc:sim}
\( \models \IT{Sim}(S^m(\Za),S^p(\Zb)) \) iff
$m = p$.
\end{LEMMA}

\begin{proof}
We have \( \models \IT{Sim}(S^m(\Za),S^p(\Zb)) \) iff
$\models \Za \doteq \Zb \RW S^m(\Za) \doteq S^p(\Zb)$
iff, by \XL{pc:cl}(a),
\( \models S^m(\Zb) \doteq S^p(\Zb) \)
iff, by \XL{pc:cl}(b),
\( S^m(\Zb) = S^p(\Zb) \) iff
$m = p$.
\end{proof}

\PAR{Addition}\label{par:arpc:add:df}
Denote by $\IT{Plus}(x,y,z)$ the semiformula
$\Zb \doteq x \RW z \doteq y$
and by $\IT{Add}(x,y,z,w)$ the semiformula
\[
\widetilde{\IT{Num}}(w) \wedge \IT{Sim}(y,w) \wedge \IT{Plus}(x,w,z) \, .
\]

\begin{LEMMA}{}\label{cl:arpc:add1}\label{par:arpc:add1}
$\models \IT{Plus}(S^m(\Za),S^p(\Zb),S^q(\Za))$
iff
$q = m+p$.
\end{LEMMA}

\begin{proof}
$\models \IT{Plus}(S^m(\Za),S^p(\Zb),S^q(\Za))$ iff
$\models \Zb \doteq S^m(\Za) \RW S^q(\Za) \doteq S^p(\Zb)$
iff, by \XL{pc:cl}(a),
\( \models S^q(\Za) \doteq S^p(S^m(\Za)) \) iff, by \XL{pc:cl}(b),
\( S^q(\Za) = S^p(S^m(\Za)) \) iff
$q = m+p$.
\end{proof}

\PAR{Proof of \XL{arpc}(b)}\label{par:arpc:add2}
\( \IT{Add}(S^m(\Za),S^p(\Za),S^q(\Za),\BO{\ast}) \) is solvable iff
\( \models \IT{Add}(S^m(\Za),S^p(\Za),S^q(\Za),\BO{d}) \)
for some \BO{d} iff
\[ \models \widetilde{\IT{Num}}(\BO{d}) \wedge \IT{Sim}(S^p(\Za),\BO{d})
             \wedge \IT{Plus}(S^m(\Za),\BO{d},S^q(\Za)) \]
for some \BO{d} iff, by \XL{arpc:a}(b),
\[ \models \IT{Sim}(S^p(\Za),S^{p_1}(\Zb))
             \wedge \IT{Plus}(S^m(\Za),S^{p_1}(\Zb),S^q(\Za)) \]
for some $p_1$ iff, by \XL{arpc:sim},
$\models \IT{Plus}(S^m(\Za),S^p(\Zb),S^q(\Za))$
iff, by \XL{arpc:add1}, $q = m+p$.
\qed

\PAR{Tables}\label{par:arpc:tab:df}
Semiterms $\BO{a}(x,y,z)$ of the form
\[
(S^{p_1}(x) \PP S^{q_1}(y)) \PP (S^{p_2}(x) \PP S^{q_2}(y)) \PP \ldots \PP
(S^{p_r}(x) \PP S^{q_r}(y)) \PP z
\]
are {\em semitables\/} of length $r \geq 0$.
Note that the term $z$ is a semitable of length $0$.
Closed instances of
semitables are {\em tables\/}.  Denote by $\IT{Tab}(x)$ the semiformula
\[
\Za \doteq S(\Za) \wedge k \doteq (\Za \PP \Za) \PP k \RW k \doteq x
\]
and by $\widetilde{\IT{Tab}}(x)$ the semiformula
\[
\Zc \doteq S(\Zc) \wedge \Zb \doteq S(\Zb) \wedge
\tilde{k} \doteq (\Zc \PP \Zb) \PP \tilde{k} \RW \tilde{k} \doteq x  \, .
\]

\begin{LEMMA}{}\label{cl:arpc:tab}\label{par:arpc:tab}
\begin{enumerate}
\item[\rm (a)]
$\IT{Tab}(\BO{\ast})$ is solved exactly by
$\BO{a}(\Za,\Za,k)$ where $\BO{a}(x,y,z)$ is a semitable,
\item[\rm (b)]
$\widetilde{\IT{Tab}}(\BO{\ast})$ is solved exactly by
$\BO{a}(\Zc,\Zb,\tilde k)$ where $\BO{a}(x,y,z)$ is a semitable.
\end{enumerate}
\end{LEMMA}
\begin{proof}
We prove only the part (b) as the proof of (a) is similar and even simpler.
$\models \widetilde{\IT{Tab}}(\BO{a}(\Zc,\Zb,\tilde k))$ is proved by a
straightforward induction on the length of
the semitable $\BO{a}(x,y,z)$.
Conversely, if
$\models \IT{Tab}(\BO{b})$ for a term \BO{b} then also
$\SM \models \IT{Tab}(\BO{b})$ for a structure
$\SM$ with the
domain $\{0,2,3,4,5\}$ and the interpretation of function symbols $I$
such that
$\Zc^I = 2$, $\Zb^I = 3$, $\tilde{k}^I = 4$,
$S^I(2) = 2$, $S^I(3) = 3$, and $S^I(d) = 0$
for all other elements $d$ of the domain.
Furthermore,
$2\PP^I 3 = 5$; $5 \PP^I 4 = 4$; and $d_1 \PP^I d_2 = 0$
for all other elements $d_1$, $d_2$ of the domain.
All other function symbols are interpreted as functions always yielding $0$.
Now, for the antecedent of $\widetilde{\IT{Tab}}(\BO{b})$ we have
\[
\SM \models \Zc \doteq S(\Zc) \wedge \Zb \doteq S(\Zb) \wedge
\tilde{k} \doteq (\Zc \PP \Zb) \PP \tilde{k}
\]
and so
$\SM \models \tilde{k} \doteq \BO{b}$, i.e. $4 = \tilde{k}^I = \BO{b}^I$.
This is clearly possible only if
the term \BO{b} is of a form
\( (\BO{b}_1 \PP \BO{c}_1) \PP (\BO{b}_2 \PP \BO{c}_2) \PP \ldots \PP
   (\BO{b}_n \PP \BO{c}_n) \PP \tilde{k} \)
for $n \geq 0$, $\BO{b}_j^I = 2$, and $\BO{c}_j^I = 3$
for all $1 \leq j \leq n$.
By a reasoning
similar to that in the proof of \XL{arpc:a}
we get that $\BO{b}_j$ and
$\BO{c}_j$ are $\Zc$ and $\Zb$-numerals respectively. Hence,
$\BO{b} = \BO{a}(\Zc,\Zb,\tilde{k})$ for a semitable $\BO{a}(x,y,z)$.
\end{proof}

\PAR{Similar tables}\label{par:arpc:simt:df}
Denote by $\widetilde{\IT{Sim}}(x,y)$ the semiformula
\[
\Za \doteq \Zc \wedge \Za \doteq \Zb \wedge k \doteq \tilde{k} \RW x \doteq y
\, .
\]

\begin{LEMMA}{}\label{cl:arpc:simt}\label{par:arpc:simt}
For semitables $\BO{a}(x,y,z)$ and $\BO{b}(x,y,z)$ we have
\[
\models \widetilde{\IT{Sim}}(\BO{a}(\Za,\Za,k),\BO{b}(\Zc,\Zb,\tilde{k}))\ 
\mbox{iff}\ 
\BO{a}(x,y,z) = \BO{b}(x,y,z) \, .
\]
\end{LEMMA}

\begin{proof}
We have
$\models \widetilde{\IT{Sim}}(\BO{a}(\Za,\Za,k),\BO{b}(\Zc,\Zb,\tilde{k}))$
iff
\[
\models
\Za \doteq \Zc \wedge \Za \doteq \Zb \wedge k \doteq \tilde{k} \RW
\BO{a}(\Za,\Za,k) \doteq \BO{b}(\Zc,\Zb,\tilde{k})
\]
iff, by \XL{pc:cl},
\( \BO{a}(\Za,\Za,k) = \BO{b}(\Za,\Za,k) \)
iff, by a straightforward induction on the length of the semitable
$\BO{b}(x,y,z)$, we have $\BO{a}(x,y,z) = \BO{b}(x,y,z)$.
\end{proof}

\PAR{\texorpdfstring{$(m,p)$}{(m,p)}-tables}\label{par:arpc:mp:df}
The tables we are interested in encode `course of values' of multiplication
and we call them {\em $(m,p)$-tables}.
The term
$\BO{a}(x,y,z)$ semitable is a {\em $(m,0)$-semitable\/} if
\( \BO{a}(x,y,z) = z \).
If $\BO{a}(x,y,z)$ is a $(m,p)$-semitable then the term
\( (S^p(x) \PP S^{m\cdot p}(y)) \PP \BO{a}(x,y,z) \)
is a {\em $(m,p+1)$-semitable\/}.
Note that $(m,p)$-semitables are of length $p$ and they differ only in the
variables $x$, $y$, and $z$.
Closed instances of $(m,p)$-semitables are {\em $(m,p)$-tables\/}.

\begin{LEMMA}{}\label{cl:arpc:mp}\label{par:arpc:mp}
For a semitable $\BO{a}(x,y,z)$ we have
\begin{equation}\label{eq:arpc:mp}
\BO{a}(S(\Za),S^m(\Za),((\Za \PP \Za) \PP k)) =
(S^p(\Za) \PP S^q(\Za)) \PP \BO{a}(\Za,\Za,k)
\end{equation}
iff $q = m \cdot p$ and $\BO{a}(x,y,z)$ is a $(m,p)$-semitable.
\end{LEMMA}
\begin{proof}
By induction on the length of the semitable $\BO{a}(x,y,z)$.
If $\BO{a}(x,y,z) = z$ then \XELerr{arpc:mp}{} holds iff
$(\Za \PP \Za) \PP k = (S^p(\Za) \PP S^q(\Za)) \PP k$
iff $p=0$ and $q = 0$ iff
$q = m \cdot p$ and $z$ is a
$(m,p)$-semitable.

If
\( \BO{a}(x,y,z) = (S^{p_1}(x) \PP S^{q_1}(y)) \PP \BO{a}_1(x,y,z) \)
for some semitable $\BO{a}_1(x,y,z)$
then
\XELerr{arpc:mp}{} holds iff, after some simplifications, we have
\[
\BO{a}_1(S(\Za),S^m(\Za),((\Za \PP \Za) \PP k)) =
(S^{p_1}(\Za) \PP S^{q_1}(\Za)) \PP \BO{a}_1(\Za,\Za,k)\, ,
\]
and
$S^{p_1+1}(\Za) \PP S^{m+q_1}(\Za) =
S^{p}(\Za) \PP S^{q}(\Za)$
iff, by IH,
\begin{equation}\label{eq:arpc:mp:1}
\mbox{
$q_1 = m \cdot p_1$, $\BO{a}_1(x,y,z)$ is a $(m,p_1)$-semitable,
and $p_1+1 = p$, $m+q_1 = q$.
}
\end{equation}
From \XEL{arpc:mp}{1} we have that
$\BO{a}(x,y,z)$ is a $(m,p_1+1)$-semitable, i.e. $(m,p)$-semitable,
and
\[ q = m + q_1 =  m + m \cdot p_1 = m \cdot (p_1+1) = m \cdot p \, . \]
Vice versa,
if $q = m \cdot p$ and
$\BO{a}(x,y,z)$ is a $(m,p)$-semitable then from the last
it must be the case that
$p = p_1+1$, $q_1 = m \cdot p_1$, and that
$\BO{a}_1$ is a $(m,p_1)$-semitable.
Moreover,
\[ q = m \cdot p = m \cdot (p_1+1) = m \cdot p_1 + m = q_1 + m  \]
and \XEL{arpc:mp}{1} holds.
\end{proof}

\PAR{Multiplication}\label{par:arpc:mul:df}
Denote by $\IT{Tim}(x,y,z,w,\tilde{w})$ the semiformula
\[
\Zc \doteq S(\Za) \wedge \Zb \doteq x \wedge \tilde{k} \doteq (\Za \PP \Za) \PP
k \RW \tilde{w} \doteq (y \PP z) \PP w
\]
and by $\IT{Mul}(x,y,z,w,\tilde{w})$ the semiformula
\[
\IT{Tab}(w) \wedge \widetilde{\IT{Tab}}(\tilde{w}) \wedge
\IT{Sim}(w,\tilde{w}) \wedge
\IT{Tim}(x,y,z,w,\tilde{w}) \, .
\]

\begin{LEMMA}{}\label{cl:arpc:mul1}\label{par:arpc:mul1}
For a semitable $\BO{a}(x,y,z)$ we have
\[ \models \IT{Tim}(S^m(\Za),S^p(\Za),S^q(\Za),
   \BO{a}(\Za,\Za,k),\BO{a}(\Zc,\Zb,\tilde{k})) \]
iff $q = m \cdot p$ and $\BO{a}(x,y,z)$ is a $(m,p)$-semitable.
\end{LEMMA}

\begin{proof}
We have
$\models \IT{Tim}(S^m(\Za),S^p(\Za),S^q(\Za),
   \BO{a}(\Za,\Za,k),\BO{a}(\Zc,\Zb,\tilde{k}))$
iff
\[
\models \Zc \doteq S(\Za) \wedge \Zb \doteq S^m(\Za) \wedge
\tilde{k} \doteq (\Za \PP \Za) \PP k \RW
\BO{a}(\Zc,\Zb,\tilde{k}) \doteq
(S^p(\Za) \PP S^q(\Za)) \PP \BO{a}(\Za,\Za,k)
\]
iff, by \XL{pc:cl}, we have \XEGerr{arpc:mp}{} iff, by \XL{arpc:mp},
$q = m \cdot p$ and $\BO{a}(x,y,z)$ is a $(m,p)$-semitable.
\end{proof}

\PAR{Proof of \XL{arpc}(c)}\label{par:arpc:mul2}
We wish to prove that
\[ \IT{Mul}(S^m(\Za),S^p(\Za),S^q(\Za),\BO{\ast},\BO{\ast}_1) \]
is solvable iff
\( \SN \models S^m(\Za) \cdot S^p(\Za) \doteq S^q(\Za) \).
The proof is similar to that of \XL{arpc}(b) as given in
\XP{arpc:add2}.
\qed

\PAR{Remark}\label{par:arpc:rem}
Note that the formula
\( \IT{Add}(S^m(\Za),S^p(\Za),S^{m+p}(\Za),\BO{\ast}) \)
has a unique solution $S^p(\Zb)$ which has the same size as $S^p(\Za)$
where under the {\em size\/} of a term we understand the number of its function
symbols.

Note also that the formula
\(
\IT{Mul}(S^m(\Za),S^p(\Za),S^{m \cdot p}(\Za),\BO{\ast},\BO{\ast}_1)
\)
has a unique solution where both $(m,p)$-tables have the same size given by a
polynomial in $m$ and $p$.


\section{Non-recursiveness of \texorpdfstring{$Sk_n(\psi)$}{Skₙ(Ψ)}}\label{sc:skn}

The reader has surely noted that for every PC-arithmetic semiformula
\( \phi(\bar{\BO{x}},\bar{\BO{w}}) \)
associated to a diophantine semiformula
\( \psi(\bar{\BO{x}}) \)
the formulas
\( \phi(\bar{\BO{\ast}},\bar{\BO{\ast}}_1) \)
are instances of SREU.  Hence, the undecidability of $\IT{Sk}_1(\psi)$ (see
\XT{sk1:col}) can be proved also from the proof of the undecidability of SREU
by Degtyarev and Voronkov \cite{qq:degvor:rigid:1}.  Actually, in a recent
paper \cite{qq:degvor:rigid:3}, which cites an earlier version of the present
paper, they gave an alternative proof of the undecidability of SREU by the
reduction of Matiyasevich's result in a similar way as we did in the previous
sections.

The undecidability of $\IT{Sk}_n(\psi)$ for $n > 1$, which is our main result
(\XT{skn:col}), requires additional work.  We now give an example illustrating
why the generalization to $\IT{Sk}_n(\psi)$ is not straightforward.

\PAR{Example}\label{par:skn:exam}
Consider the semiformula $\phi(x,w) \in \Apc$ of the form
\[
\IT{Num}(x) \wedge
\IT{Num}(S(\Za)) \wedge
\IT{Num}(\Za) \wedge
\IT{Add}(x,S(\Za),\Za,w)
\]
which is associated to the diophantine semiformula
\( x + S(\Za) \doteq \Za \).
We have $\SN \not \models \exists x\, x + S(\Za) \doteq \Za$ and so
the formula
$\phi(\ast,\ast_1)$,
which is a 1-skeleton of
\( \exists x \exists w\, \phi(x,w) \),
is not solvable by \XL{sk1:cl3}.
On the other hand, we can solve a $2$-skeleton of the last formula
because we have
\begin{equation}\label{eq:skn:exam:2}
\models \phi(\Za,\Zb) \vee \phi(S(\Za),S(\Zb))\, .
\end{equation}
Indeed, \XEL{skn:exam}{2} holds iff,
after some simplifications which remove valid
subformulas,
\begin{equation}\label{eq:skn:exam:3}
\models \IT{Sim}(S(\Za),\Zb) \vee \IT{Plus}(S(\Za),S(\Zb),\Za)
\end{equation}
iff
\[
\models
(\Za \doteq \Zb \RW S(\Za) \doteq \Zb) \vee
(\Zb \doteq S(\Za) \RW \Za \doteq S(\Zb))
\]
iff
\begin{equation}\label{eq:skn:exam:4}
\models
\Za \doteq \Zb \wedge \Zb \doteq S(\Za) \RW
S(\Za) \doteq \Zb \vee \Za \doteq S(\Zb)
\end{equation}
and the last is the case.

This is a typical situation where the solvability of a disjunction does not
guarantee the solvability of disjuncts because their clauses may interfere
through solutions (compare \XEL{skn:exam}{3} to \XEL{skn:exam}{4}) even if
the disjuncts do not share the unknowns.

\PAR{Variants of PC-arithmetic semiformulas}\label{par:skn:Ai}
For each $i \geq 1$
denote by $P_i$
the language consisting
of function symbols $S$, `$\PP$'
and of constants
$\Za_i$, $\Zc_i$, $\Zb_i$, $k_i$, $\tilde{k}_i$.
The constants of languages $P_i$ are called {\em
special\/} constants.
For simplicity sake we use a different set of variables
for each language $P_i$:
\( x^i_1, x^i_2, x^i_3, \ldots .\)
We
use $\BO{x}_i$ and $\bar{\BO{x}}_i$ as metavariables ranging over
variables and sequences of variables of $P_i$ respectively.
We do not distinguish between the numeric and
table variables of PC-arithmetic semiformulas which will be written in the form
$\phi(\bar{\BO{x}})$.
The semiformula
\( \phi_i(\bar{\BO{x}}_i) \)
of the language $P_i$ is called a {\em variant\/} of a PC-arithmetic
semiformula
\( \phi(\bar{\BO{x}})\)
if it is obtained from the last semiformula by replacing each of its
constants \BO{k} by the special constant $\BO{k}_i$ and
each of its variables $x_j$ by $x^i_j$.
Clearly, $\phi(\bar{\BO{\ast}})$ is solvable iff $\phi_i(\bar{\BO{\ast}})$ is.
Moreover, the solutions differ only in the corresponding constants. We denote
by $\Apc_i$ the class of semiformulas of $P_i$ which are variants of
semiformulas of $\Apc$.

For $n \geq 1$ we assign to every quantifier-free semiformula
\( \phi(\bar{\BO{x}}) \in \Apc \)
the semiformula
\( \psi(\bar{\BO{x}}_1,\ldots,\bar{\BO{x}}_n) \)
of the form
\[
\phi_1(\bar{\BO{x}}_1) \wedge \cdots \wedge \phi_n(\bar{\BO{x}}_n)\, ,
\]
where every
\( \phi_i(\bar{\BO{x}}_i) \in \Apc_i \)
for $1 \leq i \leq n$
is a variant of $\phi(\bar{\BO{x}})$.
We denote by $\Apc^{(n)}$ the class of
semiformulas obtained by this assignment from the class $\Apc$.

\PAR{Invariancy of assignment under substitution}\label{par:skn:sub}
We will use the following fact which is similar to that in \XP{sk1:sub}:
\begin{quote}
if for $n \geq 1$ the semiformula
\(
\psi(\BO{x}_1,\bar{\BO{x}}_1,\ldots,\BO{x}_n,\bar{\BO{x}}_n) \in \Apc^{(n)}
\)
is assigned to the semiformula
\( \phi(\BO{x},\bar{\BO{x}}) \in \Apc \)
where \BO{x} is a numeric variable then for every $\Za$-numeral $S^m(\Za)$
the semiformula
\(
\psi(S^m(\Za_1),\bar{\BO{x}}_1,\ldots,S^m(\Za_n),\bar{\BO{x}}_n) \in \Apc^{(n)}
\)
is assigned to the semiformula
\( \phi(S^m(\Za),\bar{\BO{x}}) \in \Apc \)
\end{quote}
which is easily proved by induction on
\( \phi(\BO{x},\bar{\BO{x}}) \).

The Main theorem \ref{cl:skn:col} on the undecidability of $n$-skeletons
requires two lemmas the first of which will be proved in the next section.

\PAR{Main lemma}\label{cl:skn:mainlemma}\label{par:skn:mainlemma}
{\em
If for $n \geq 1$ the semiformulas
\( \phi_i(\bar{\BO{x}}_i) \in \Apc_i \)
are variants of
\( \phi(\bar{\BO{x}}) \in \Apc \)
for all
\( 1 \leq i \leq n \)
and
\begin{equation}\label{eq:skn:mainlemma:1}
\models \phi_1(\bar{\BO{a}}_1) \vee \cdots \vee \phi_n(\bar{\BO{a}}_n)
\end{equation}
then
\( \models \phi_i(\bar{\BO{a}}_i) \)
for some $i$.
}

\paragraph*{Proof outline.}
If
\( \not \models \phi_i(\bar{\BO{a}}_i) \)
for all $1 \leq i \leq n$ then we will construct a structure $\SM$ falsifying
the variants at the same time:
\( \SM \not \models \phi_i(\bar{\BO{a}}_i) \)
for all $1 \leq i \leq n$. We will then have
\(
\SM \not \models \phi_1(\bar{\BO{a}}_1) \vee \cdots \vee \phi_n(\bar{\BO{a}}_n)
\)
and \XEL{skn:mainlemma}{1} will not hold.  We postpone the construction of
$\SM$ until \XP{main:proof}.
\qed

\begin{LEMMA}{}\label{cl:skn:sk1}\label{par:skn:sk1}
If for $n \geq 1$ the semiformula
\( \psi(\bar{\BO{x}}_1,\ldots,\bar{\BO{x}}_n) \in \Apc^{(n)} \)
is assigned to \( \phi(\bar{\BO{x}}) \in \Apc \)
then
\[
\IT{Sk}_1(\exists \bar{\BO{x}} \phi(\bar{\BO{x}}))\
\mbox{iff}\
\IT{Sk}_n(\exists \bar{\BO{x}}_1 \ldots \exists \bar{\BO{x}}_n\,
                        \psi(\bar{\BO{x}}_1,\ldots,\bar{\BO{x}}_n))\, .
\]
\end{LEMMA}

\begin{proof}
($\mRW$):
If
\( \IT{Sk}_1(\exists \bar{\BO{x}}\, \phi(\bar{\BO{x}})) \)
then $\phi(\bar{\BO{\ast}})$ is
solved by some terms $\bar{\BO{a}}$ of $P$ and
we can solve
\( \psi(\bar{\BO{\ast}}_1,\ldots,\bar{\BO{\ast}}_n) \)
by the corresponding variants $\bar{\BO{a}}_i$.
Thus also any $n$-skeleton of
\(
\exists \bar{\BO{x}}_1 \ldots \exists \bar{\BO{x}}_n\,
\psi(\bar{\BO{x}}_1,\ldots,\bar{\BO{x}}_n)
\)
is solvable.

($\mLW$):
If some $n$-skeleton of
\(
\exists \bar{\BO{x}}_1 \ldots \exists \bar{\BO{x}}_n\,
\psi(\bar{\BO{x}}_1,\ldots,\bar{\BO{x}}_n)
\)
is solvable then
\[
\models \psi_1(\bar{\BO{a}}_{1,1},\ldots,\bar{\BO{a}}_{1,n}) \vee \cdots \vee
                        \psi_n(\bar{\BO{a}}_{n,1},\ldots,\bar{\BO{a}}_{n,n})
\]
for some terms $\bar{\BO{a}}_{j,k}$ ($1 \leq j,k \leq n$). Distributing
$\vee$'s over $\wedge$'s and weakening yields
\[
\models \phi_1(\bar{\BO{a}}_{1,1}) \vee \cdots \vee
                        \phi_n(\bar{\BO{a}}_{n,n}) \, .
\]
By the Main lemma \ref{cl:skn:mainlemma} we have
\( \models \phi_j(\bar{\BO{a}}_{j,j}) \)
for some $j$. Hence, $\phi_j(\bar{\BO{\ast}}_j)$
and also $\phi(\bar{\BO{\ast}})$
are solvable.
\end{proof}

\begin{THEOREM}{}\label{cl:skn:cl}\label{par:skn:cl}
Let $n \geq 1$. To every recursively enumerable predicate $R(m)$ there is a
semiformula
\(
\psi(\BO{x}_1,\bar{\BO{x}}_1,\ldots,\BO{x}_n,\bar{\BO{x}}_n) \in \Apc^{(n)}
\)
such that for all $m$
\begin{equation}\label{eq:skn:cl:1}
R(m)\ \mbox{iff}\
\IT{Sk}_n(\exists \bar{\BO{x}}_1 \ldots \exists \bar{\BO{x}}_n\,
\psi(S^m(\Za_1),\bar{\BO{x}}_1,\ldots,S^m(\Za_n),\bar{\BO{x}}_n))\, .
\end{equation}
\end{THEOREM}

\begin{proof}
By \XT{sk1:cl} there is a PC-arithmetic semiformula
\( \phi(\BO{x},\bar{\BO{x}}) \in \Apc \)
such that for every number $m$
\begin{equation}\label{eq:skn:cl:2}
R(m)\ \mbox{iff}\
\IT{Sk}_1(\exists \bar{\BO{x}}\, \phi(S^m(\Za),\bar{\BO{x}})) \, .
\end{equation}
Take the semiformula
\(
\psi(\BO{x}_1,\bar{\BO{x}}_1,\ldots,\BO{x}_n,\bar{\BO{x}}_n) \in \Apc^{(n)}
\)
assigned to
\( \phi(\BO{x},\bar{\BO{x}}) \in \Apc \).
Then $R(m)$ holds iff, by \XEL{skn:cl}{2},
\( \IT{Sk}_1(\exists \bar{\BO{x}}\,\phi(S^m(\Za),\bar{\BO{x}})) \)
iff, by \XP{skn:sub} and \XL{skn:sk1},
\(
\IT{Sk}_n(\exists \bar{\BO{x}}_1 \ldots \exists \bar{\BO{x}}_n\,
\psi(S^m(\Za_1),\bar{\BO{x}}_1,\ldots,S^m(\Za_n),\bar{\BO{x}}_n))
\).
\end{proof}

\PAR{Main theorem}{}\label{cl:skn:col}\label{par:skn:col}
{\em
The predicate $\IT{Sk_n}(\psi)$ is not recursive for any $n \geq 1$.
}

\begin{proof}
We fix $n$, take any recursively enumerable but not recursive predicate $R(m)$,
and obtain a semiformula
\(
\psi(\BO{x}_1,\bar{\BO{x}}_1,\ldots,\BO{x}_n,\bar{\BO{x}}_n) \in \Apc^{(n)}
\)
from \XT{skn:cl}. We can clearly find a primitive recursive function $f$ such
that
\[
f(m) =
\exists \bar{\BO{x}}_1 \ldots \exists \bar{\BO{x}}_n\,
\psi(S^m(\Za_1),\bar{\BO{x}}_1,\ldots,S^m(\Za_n),\bar{\BO{x}}_n)\, .
\]
Then $R(m)$ iff $\IT{Sk}_n(f(m))$ by \XEG{skn:cl}{1}. If the predicate
$\IT{Sk}_n(\psi)$ were recursive so would be $R(m)$.
\end{proof}

\PAR{Discussion}\label{par:skn:disc}
From \XT{skn:col} we can immediately see that for no $n$ there can be a
recursive function $f_n(\psi)$ which yields a bound on the size of solutions to
$n$-skeletons of an existential formula $\psi$. Existence of such a function would
make the non-recursive predicate $\IT{Sk}_n(\psi)$ recursive by a simple test
of all candidate solutions of the size less than the bound.


\section{Proof of the Main Lemma}\label{sc:main}

For the proof of the Main lemma \ref{cl:skn:mainlemma} given in
\XP{main:proof} we must be able to falsify simultaneously $n$ variants
of a not valid formula simulating arithmetic in predicate calculus.
For this we need a structure interpreting the special constants
in such a way that the interpretations do not interfere with each other.
Such structures form a family defined in the following paragraph.

\PAR{Structures \texorpdfstring{$\MN$}{𝓜α}}\label{par:main:struct}
We define a family of structures $\MN$ with natural numbers as domains.
The structures are parameterized by functions
$\alpha(\BO{k})$ interpreting the special constants into natural numbers.

We partition the domain of natural numbers into seven mutually disjoint
subsets.  This is done with the help of a pairing
function
\[ J(j,k) = (j+k)(j+k+1)+k+1 \]
which is the standard pairing function as used in recursion theory (see for
instance \cite[pg. 43]{qq:davis:58}) but offset by one.  Thus $0$ is the only
number not in the range of $J$ and we have the {\em pairing\/} property that
from $J(i,k) = J(i',k')$ we get $i = i'$ and $k = k'$.
A computer scientist will realize that $J$ can ve viewed as the
function $\IT{cons}$ of LISP, $0$ as \IT{nil}, and the set of natural
numbers as the set of S-expressions generated from the single atom \IT{nil}.

The first set in the partition is the set
\( \{ J(0,i) \mid i \in N \} \)
and it plays the role of natural numbers, where the natural number $i$ is
embedded into $J(0,i)$.  The next five sets are
\( \{ J(j,i) \mid i \in N \} \)
for $1 \leq j \leq 5$ and they will be used to interpret the special constants
in a special way.  The seventh set consists of the remaining natural numbers,
i.e. of $0$ and of the numbers of the form $J(j+6,i)$, and it will play no
special role.

For a given function $\alpha$
we define the interpretation $I$ of the structure $\MN$
as follows. All predicate symbols \BO{p} are
always false: $\BO{p}^I = \emptyset$.
For a special
constant $\BO{k}$ we have $\BO{k}^I = \alpha(\BO{k})$.  Function
symbols \BO{f}
other than the special constants, $S$, and `$\PP$' are interpreted as
functions always yielding $0$: $\BO{f}^I(\bar{d}) = 0$.
The interpretation of the unary function symbol $S$ is as follows:
\begin{eqnarray*}
S^I\,J(0,i) & = & J(0,i+1)
\\
S^I\,J(1,i) & = & J(1,i)
\\
S^I\,J(2,i) & = & J(2,i)
\\
S^I\,J(3,i) & = & J(3,i)
\\
S^I(m) & = & 0 \ \mbox{otherwise.}
\end{eqnarray*}
We can see that the function $S^I$ behaves as the successor function on the
numbers of the form $J(0,i)$, while on the other numbers it has special
properties.

The interpretation of the binary function symbol `$\PP$' is as follows:
\begin{eqnarray*}
J(0,i) \PP^I J(0,j) & = & J(0,J(i,j))
\\
J(1,i) \PP^I J(1,i) & = & J(5,i)
\\
J(2,i) \PP^I J(3,i) & = & J(5,i)
\\
J(5,i) \PP^I J(4,i) & = & J(4,i)
\\
m \PP^I n & = & 0 \ \mbox{otherwise.}
\end{eqnarray*}
We can see that the function `$\PP^I$' {\em reflects\/} $J$, i.e. behaves as a pairing
function, on the numbers of the form $J(0,i)$.

The following simple property of structures will be
often used in the proofs below.

\PAR{Lemma}\label{cl:main:sub}\label{par:main:sub}
{\em
Let $\SM$ be a structure and let $\phi(\BO{x})$ be a semiformula with
at most $\BO{x}$ free. If the terms
$\BO{a}$ and
$\BO{b}$ denote in $\SM$ the same element, i.e. if
$\SM \models \BO{a} \doteq \BO{b}$, then
\( \SM \models \phi(\BO{a}) \LRW \phi(\BO{b}) \).
}

\begin{proof}
See the lemma on Pg.~19 in \cite{qq:shoenfield}.
\end{proof}

\PAR{Note}\label{par:main:variants}
By the conventions on variants
of PC-arithmetic semiformulas discussed in \XP{skn:Ai} we denote
by $\IT{Num}_i(x) \in \Apc_i$ the corresponding variant of
$\IT{Num}(x) \in \Apc$. Similarly, for other semiformulas defined
in \XS{arpc}. Lemmas proved in \XS{arpc} hold also for
the corresponding variants. To emphasize this we will, for instance,
refer to
\XL{arpc:a}$_i$ when we mean
\XL{arpc:a} modified for the variant $\IT{Num}_i(x)$.

\begin{LEMMA}{}\label{cl:main:num}\label{par:main:num}
\begin{enumerate}
\item[\rm (a)]
\( \models \IT{Num}_i(\BO{a}) \) iff
\( \MN \models \IT{Num}_i(\BO{a}) \)
where
\( \alpha(\Za_i) = J(1,i) \)
and
\( \alpha(\BO{k}) \neq J(1,i) \)
for all other special constants,

\item[\rm (b)]
\( \models \widetilde{\IT{Num}}_i(\BO{a}) \) iff
\( \MN \models \widetilde{\IT{Num}}_i(\BO{a}) \)
where
\( \alpha(\Zb_i) = J(3,i) \)
and
\( \alpha(\BO{k}) \neq J(3,i) \)
for all other special constants.
\end{enumerate}
\end{LEMMA}

\begin{proof}
We prove only the part (a) as the proof of (b) is similar.
The direction $\mRW$ is obvious.
For the direction $\mLW$ we assume
\( \not \models \IT{Num}_i(\BO{a}) \) and then
\( \BO{a} = S^m\,\BO{f}(\bar{\BO{b}}) \)
for some $m$, $\bar{\BO{b}}$, and $\BO{f}$ which is neither $S$ nor $\Za_i$
by \XL{arpc:a}(a)$_i$.
We can see similarly as in the proof of \XL{arpc:a}(a)
that \( \MN \models \Za_i \doteq S(\Za_i) \) and
\( \MN \not\models \Za_i \doteq S^{m}\,\BO{f}(\bar{\BO{b}}) \).
Hence, \( \MN \not \models \IT{Num}_i(\BO{a}) \).
\end{proof}

\begin{LEMMA}{}\label{cl:main:sim}\label{par:main:sim}
\( \models \IT{Sim}_i(S^m(\Za_i),S^p(\Zb_i)) \) iff
\( \MN \models \IT{Sim}_i(S^m(\Za_i),S^p(\Zb_i)) \)
where
\( \alpha(\Za_i) = \alpha(\Zb_i) = J(0,0) \).
\end{LEMMA}

\begin{proof}
The direction $\mRW$ is obvious.
For the direction $\mLW$ we assume
\( \MN \models \IT{Sim}_i(S^m(\Za_i),S^p(\Zb_i)) \) and, since
\( \Za_i^I = J(0,0) = \Zb_i^I \),
we have
\( \MN \models S^m(\Zb_i) \doteq S^p(\Zb_i) \)
by \XL{main:sub}.
$S^I$ behaves as the
successor function on the numbers of the form $J(0,j)$ and thus
it must be the case that $m = p$. We now get
\( \models \IT{Sim}_i(S^m(\Za_i),S^p(\Zb_i)) \)
by \XL{arpc:sim}$_i$.
\end{proof}

\begin{LEMMA}{}\label{cl:main:add1}\label{par:main:add1}
\( \models \IT{Plus}_i(S^m(\Za_i),S^p(\Zb_i),S^q(\Za_i)) \) iff
\[ \MN \models \IT{Plus}_i(S^m(\Za_i),S^p(\Zb_i),S^q(\Za_i)) \]
where
\( \alpha(\Za_i) = J(0,0) \)
and
\( \alpha(\Zb_i) = J(0,m) \).
\end{LEMMA}

\begin{proof}
The direction $\mRW$ is obvious.
For the direction $\mLW$ we assume
\( \MN \models \IT{Plus}_i(S^m(\Za_i),S^p(\Zb_i),S^q(\Za_i)) \),
i.e.
\[
\MN \models \Zb_i \doteq S^m(\Za_i) \RW S^q(\Za_i) \doteq S^p(\Zb_i)
\, .
\]
We can easily see that
$(S^m(\Za_i))^I = J(0,m) = \Zb_i^I$
and hence
\( \MN \models S^q(\Za_i) \doteq S^p(S^m(\Za_i)) \)
by \XL{main:sub}. By similar arguments as in the proof of \XL{main:sim} we can
see that $q=m+p$. We now get
\( \models \IT{Plus}_i(S^m(\Za_i),S^p(\Zb_i),S^q(\Za_i)) \)
by \XL{arpc:add1}$_i$.
\end{proof}

\begin{LEMMA}{}\label{cl:main:tab}\label{par:main:tab}
\begin{itemize}
\item[\rm (a)]
\( \models \IT{Tab}_i(\BO{a}) \) iff
\( \MN \models \IT{Tab}_i(\BO{a}) \)
where
\( \alpha(\Za_i) = J(1,i) \),
\( \alpha(k_i) = J(4,i) \),
and for all other special constants \BO{k} we have
\( \alpha(\BO{k}) \neq J(j+1,i) \)
for all $j$,
\item[\rm (b)]
\( \models \widetilde{\IT{Tab}}_i(\BO{a}) \) iff
\( \MN \models \widetilde{\IT{Tab}}_i(\BO{a}) \)
where
\( \alpha(\Zc_i) = J(2,i)$, $\alpha(\Zb_i) = J(3,i) \),
\( \alpha(\tilde{k}_i) = J(4,i) \),
and for all other special constants \BO{k} we have
\( \alpha(\BO{k}) \neq J(j+1,i) \)
for all $j$.
\end{itemize}
\end{LEMMA}

\begin{proof}
We prove only the part (b) as the proof of (a) is similar and even simpler.
The direction $\mRW$ is obvious.
For the direction $\mLW$ we assume
\( \MN \models \widetilde{\IT{Tab}}_i(\BO{a}) \).
Since
\[
\MN \models
\Zc_i \doteq S(\Zc_i) \wedge \Zb_i \doteq S(\Zb_i) \wedge
        \tilde{k}_i \doteq (\Zc_i,\Zb_i),\tilde{k}_i
\]
we must also have
\( \MN \models \tilde{k}_i \doteq \BO{a} \).
Reasoning similar to that
in the proof of \XL{arpc:tab}(b) shows that
\( \BO{a} = \BO{b}(\Zc_i,\Zb_i,\tilde{k}_i) \)
for a semitable $\BO{b}(x,y,z)$.
We now get
\( \models \widetilde{\IT{Tab}}_i(\BO{a}) \)
by \XL{arpc:tab}(b)$_i$.
\end{proof}

\begin{LEMMA}{}\label{cl:main:simt}\label{par:main:simt}
For semitables $\BO{a}(x,y,z)$ and $\BO{b}(x,y,z)$ we have
\begin{equation}\label{eq:main:simt:1}
\models
\widetilde{\IT{Sim}}_i(\BO{a}(\Za_i,\Za_i,k_i),\BO{b}(\Zc_i,\Zb_i,\tilde{k}_i))
\end{equation}
iff
\(
\MN \models
\widetilde{\IT{Sim}}_i(\BO{a}(\Za_i,\Za_i,k_i),\BO{b}(\Zc_i,\Zb_i,\tilde{k}_i))
\)
where
\[
\alpha(\Za_i) = \alpha(\Zc_i) = \alpha(\Zb_i) = \alpha(k_i) =
\alpha(\tilde{k}_i) = J(0,0)\, .
\]
\end{LEMMA}

\begin{proof}
The direction $\mRW$ is obvious.
For the direction $\mLW$ we assume
\(
\MN \models
\widetilde{\IT{Sim}}_i(\BO{a}(\Za_i,\Za_i,k_i),\BO{b}(\Zc_i,\Zb_i,\tilde{k}_i))
\).
Since
\[
\MN \models
\Za_i \doteq \Zc_i \wedge \Za_i \doteq \Zb_i \wedge k_i \doteq \tilde{k}_i
\]
we obtain
\(
\MN \models
\BO{a}(\Za_i,\Za_i,k_i) \doteq \BO{b}(\Za_i,\Za_i,k_i)
\)
by \XL{main:sub}. Note that all subterms of $\BO{a}(\Za_i,\Za_i,k_i)$ and
$\BO{b}(\Za_i,\Za_i,k_i)$ denote numbers of the form $J(0,j)$ on which $S$ and
`$\PP$' behave as successor and pairing functions respectively. Thus by a
straightforward induction on the semitable $\BO{b}(x,y,z)$ we get that
\( \BO{a}(x,y,z) = \BO{b}(x,y,z) \).
We now get \XEL{main:simt}{1} by \XL{arpc:simt}$_i$.
\end{proof}

\begin{LEMMA}{}\label{cl:main:mul1}\label{par:main:mul1}
For a semitable $\BO{a}(x,y,z)$ we have
\begin{equation}\label{eq:main:mul1}
\models \IT{Tim}_i(S^m(\Za_i),S^p(\Za_i),S^q(\Za_i),
        \BO{a}(\Za_i,\Za_i,k_i),\BO{a}(\Zc_i,\Zb_i,\tilde{k}_i))
\end{equation}
iff
\begin{equation}\label{eq:main:mul1:1}
 \MN \models \IT{Tim}_i(S^m(\Za_i),S^p(\Za_i),S^q(\Za_i),
        \BO{a}(\Za_i,\Za_i,k_i),\BO{a}(\Zc_i,\Zb_i,\tilde{k}_i))
\end{equation}
where
\( \alpha(\Za_i) = \alpha(k_i) = J(0,0) \),
\( \alpha(\Zc_i) = J(0,1) \),
\( \alpha(\Zb_i) = J(0,m) \),
and
$\alpha(\tilde{k}_i) = J(0,J(J(0,0),0))$.
\end{LEMMA}

\begin{proof}
The direction $\mRW$ is obvious.
For the direction $\mLW$ we assume \XEL{main:mul1}{1}, i.e.
\begin{eqnarray*}
\MN & \models &
\Zc_i \doteq S(\Za_i) \wedge \Zb_i \doteq S^m(\Za_i) \wedge
\tilde{k}_i \doteq (\Za_i \PP \Za_i) \PP k_i \RW
\\
& &
\BO{a}(\Zc_i,\Zb_i,\tilde{k}_i) \doteq
(S^p(\Za_i) \PP S^q(\Za_i)) \PP \BO{a}(\Za_i,\Za_i,k_i) \, .
\end{eqnarray*}
Since
$\Zb_i^I = J(0,1) = S(\Za_i)^I$,
$\Zc_i^I = J(0,m) = (S^m(\Za_i))^I$, and
\[ \tilde{k}_i^I = J(0,J(J(0,0),0)) =
  (\Za_i \PP \Za_i) \PP k)^I
\]
we get
\begin{equation}\label{eq:main:mul1:2}
\MN \models
\BO{a}(S(\Za_i),S^m(\Za_i),((\Za_i \PP \Za_i) \PP k_i)) \doteq
(S^p(\Za_i) \PP S^q(\Za_i)) \PP \BO{a}(\Za_i,\Za_i,k_i)
\end{equation}
by \XL{main:sub}.
Note that all subterms in \XEL{main:mul1}{2}
denote numbers of the form $J(0,j)$ on which $S$ and
`$\PP$' behave as successor and pairing functions respectively.
Thus
\[
\BO{a}(S(\Za_i),S^m(\Za_i),((\Za_i \PP \Za_i) \PP k_i)) =
(S^p(\Za_i) \PP S^q(\Za_i)) \PP \BO{a}(\Za_i,\Za_i,k_i)
\]
and we have $q = m \cdot p$ and that $\BO{a}(x,y,z)$ is a $(m,p)$-semitable
by \XL{arpc:mp}$_i$ from which we get \XELerr{main:mul1}{} by \XL{arpc:mul1}$_i$.
\end{proof}

\PAR{Proof of the Main lemma \ref{cl:skn:mainlemma}}\label{par:main:proof}
For a given $n \geq 1$ and formulas
$\phi_i(\bar{\BO{a}}_i)$
such that \( \not \models \phi_i(\bar{\BO{a}}_i) \)
for $1 \leq i \leq n$
we wish to find a structure $\MN$ falsifying all formulas:
$\MN \not \models \phi_i(\bar{\BO{a}}_i)$.
We construct such an interpretation of special
constants $\alpha$ by stages.

We set $\alpha(\BO{k}) = 0$ for all constants $\BO{k}$ of languages
$P_{n+i+1}$. The assignment of interpretations to the special constants
of $P_i$ ($1 \leq i \leq n$) is such that they never receive
interpretations $J(j+1,k)$ for any $j$ and $k \neq i$.
We now let $i$ range from $1$ through $n$. For each
$i$ exactly one of the following cases applies. The assignment $\alpha$
is constructed by stages such that in the stage $i$ the special constants
of the language $P_i$ are assigned interpretation.

{\bf (i):}
Suppose that
for a conjunct of $\phi_i(\bar{\BO{a}}_i)$ we have
\[
\not \models \IT{Num}_i(\BO{a})\
\mbox{or}\
\not \models \IT{Tab}_i(\BO{a}) \, .
\]
Setting
\( \alpha(\Za_i) = J(1,i) \),
\( \alpha(k_i) = J(4,i) \),
and
\( \alpha(\Zb_i) = \alpha(\Zc_i) = \alpha(\tilde{k}_i) = 0 \)
will cause at the end $\MN \not \models \phi_i(\bar{\BO{a}}_i)$
by Lemmas \ref{cl:main:num}(a) or \ref{cl:main:tab}(a).

{\bf (ii):}
Suppose that
{\bf (i)} does not apply and that
for a conjunct of $\phi_i(\bar{\BO{a}}_i)$ we have
\[
\not \models \widetilde{\IT{Num}}_i(\BO{a})\
\mbox{or}\
\not \models \widetilde{\IT{Tab}}_i(\BO{a}) \, .
\]
Setting
\( \alpha(\Zc_i) = J(2,i) \),
\( \alpha(\Zb_i) = J(3,i) \),
\( \alpha(\tilde{k}_i) = J(4,i) \),
and
\( \alpha(\Za_i) = \alpha(k_i) = 0 \)
will cause at the end $\MN \not \models \phi_i(\bar{\BO{a}}_i)$
by Lemmas \ref{cl:main:num}(b) or \ref{cl:main:tab}(b).

{\bf (iii):}
Suppose that
{\bf (i)} and {\bf (ii)} do not apply and that
for a conjunct of $\phi_i(\bar{\BO{a}}_i)$ we have
\[
\not \models \IT{Sim}_i(S^m(\Za_i),S^p(\Zb_i))\
\mbox{or}\
\not \models
\widetilde{\IT{Sim}}_i(\BO{a}(\Za_i,\Za_i,k_i),\BO{b}(\Zc_i,\Zb_i,\tilde{k}_i))
\, .
\]
Setting
\( \alpha(\BO{k}) = J(0,0) \)
for all special constants $\BO{k}$ of $P_i$
will cause at the end $\MN \not \models \phi_i(\bar{\BO{a}}_i)$
by Lemmas \ref{cl:main:sim} or \ref{cl:main:simt}.

{\bf (iv):}
Suppose that
{\bf (i)}, {\bf (ii)}, and {\bf (iii)} do not apply. Then
for a conjunct of $\phi_i(\bar{\BO{a}}_i)$ we must have
\begin{eqnarray*}
\mbox{\ } &
\not \models \IT{Plus}_i(S^m(\Za_i),S^p(\Za_i),S^q(\Za_i),S^p(\Zb_i))\
\mbox{or}
& \mbox{\ }
\\
\mbox{\ } &
\not \models
\IT{Tim}_i(S^m(\Za_i),S^p(\Za_i),S^q(\Za_i),
        \BO{a}(\Za_i,\Za_i,k_i),\BO{a}(\Zc_i,\Zb_i,\tilde{k}_i)) \, .
& \mbox{\ }
\end{eqnarray*}
Setting
\( \alpha(\Za_i) = \alpha(k_i) = J(0,0) \),
\( \alpha(\Zc_i) = J(0,1)$, $\alpha(\Zb_i) = J(0,m) \),
and
\( \alpha(\tilde{k}_i) = J(0,J(J(0,0),0)) \)
will cause at the end $\MN \not \models \phi_i(\bar{\BO{a}}_i)$
by Lemmas \ref{cl:main:add1} or \ref{cl:main:mul1}.
\qed


\section{Conclusion}
The negative solution of the problem of Herbrand skeletons has
important consequences for ATP.  Although ATP
involves r.e. functions we wish to use a recursive proof function which either
finds a proof or indicates reasons why a proof could not be
found. It is not very pleasant to abort a proof search
based on an r.e. function because the abort does not yield any indication
as to what has caused the failure.

Our main result shows that it is not sufficient to specify the size
of a skeleton, i.e. the number of existential axioms, one also needs
a bound on the size of terms in a solution. Thus we are interested
in finding as
efficient as possible an algorithm which, given the size of
the skeleton and a bound on the size of solutions, either finds a solution
or reports a failure. It is hoped that from the failure we can then obtain
an indication as to why the bounds were exceeded.


\newcommand{\noopsort}[1]{} \newcommand{\singletter}[1]{#1}


\begin{thebibliography}{GNPS90}

\bibitem[Bib82]{qq:bibel:82}
W.~Bibel.
\newblock {\em Automated Theorem Proving}.
\newblock Vieweg Verlag, 1982.

\bibitem[BJ80]{qq:boolos}
G.~S. Boolos and R.~C. Jeffrey.
\newblock {\em Computability and Logic}.
\newblock Cambridge University Press, second edition, 1980.

\bibitem[Bus95a]{qq:buss:intro}
S.~R. Buss.
\newblock An introduction to proof theory.
\newblock Typeset manuscript, to appear in {\em Handbook of Proof Theory\/},
  1995.

\bibitem[Bus95b]{qq:buss:herb}
S.~R. Buss.
\newblock On {H}erbrand's theorem.
\newblock Typeset manuscript, to appear in Proceedings of LCC'95, 1995.

\bibitem[Dav58]{qq:davis:58}
M.~Davis.
\newblock {\em Computability and Unsolvability}.
\newblock McGraw-Hill, 1958.

\bibitem[DV95a]{qq:degvor:rigid:1}
A.~Degtyarev and A.~Voronkov.
\newblock Simultaneous rigid {$E$}-unification is undecidable.
\newblock UPMAIL Technical Report 105, Uppsala University, Computing Science
  Department, May {\noopsort{a}}1995.

\bibitem[DV95b]{qq:degvor:rigid:3}
A.~Degtyarev and A.~Voronkov.
\newblock Simultaneous rigid {$E$}-uni\-fi\-ca\-tion is undecidable.
\newblock In {\em Proceedings of 8th Workshop, CSL'95}, {\noopsort{c}}1995.

\bibitem[Fit90]{qq:fitting}
M.~Fitting.
\newblock {\em First-order Logic and Automated Theorem Proving}.
\newblock Springer Verlag, 1990.

\bibitem[GNPS88]{qq:GNPS:88}
J.~H. Gallier, P.~Narendran, D.~Plaisted, and W.~Snyder.
\newblock Rigid {$E$}-uni\-fi\-ca\-tion is {NP}-complete.
\newblock In {\em Proceedings, Third Annual Symposium on Logic in Computer
  Science}, pages 218--227, Edinburgh, Scotland, 1988. IEEE Computer Society.

\bibitem[GNPS90]{qq:GNPS:90}
J.~H. Gallier, P.~Narendran, D.~Plaisted, and W.~Snyder.
\newblock Rigid {$E$}-uni\-fi\-ca\-tion: {NP}-completeness and applications to
  equational matins.
\newblock {\em Information and Computation}, 87(1/2):129--195, 1990.

\bibitem[Hei67]{qq:frege2godel}
J.~van Heijenoort, editor.
\newblock {\em From Frege to G{\"{o}}del: A Source Book in Mathematical Logic,
  1879-1931}.
\newblock Harvard University Press, 1967.

\bibitem[Her30]{qq:herbrand}
J.~Herbrand.
\newblock {\em Recherches sur la th{\'{e}}orie de la d{\'{e}}monstration}.
\newblock PhD thesis, University of Paris, 1930.
\newblock English translation of chapter 5 in {\cite{qq:frege2godel}}, pp.
  525--581.

\bibitem[HP93]{qq:hajekpudlak}
P.~H\'{a}jek and P.~Pudl\'{a}k.
\newblock {\em Metamathematics of First-Order Arithmetic}.
\newblock Springer Verlag, 1993.

\bibitem[Kog95]{qq:kogel}
Eric~de Kogel.
\newblock Rigid {$E$}-unification simplified.
\newblock In {\em Proceedings of 4th International Workshop, TABLEAUX'95},
  pages 17--30. Springer Verlag, 1995.

\bibitem[KP88]{qq:krajpudlak}
J.~Kraj{\'\i{}}\v{c}ek and P.~Pudl\'{a}k.
\newblock The number of proof lines and the size of proofs in first order
  logic.
\newblock {\em Arch. Math. Logic}, 27:69--84, 1988.

\bibitem[Mat70]{qq:mat}
Yu.~V. Matiyasevich.
\newblock Recursively enumerable sets are diophantine.
\newblock {\em Soviet Mathematical Doklady}, 11(2):354--358, 1970.

\bibitem[Rob65]{qq:robinson}
J.~A. Robinson.
\newblock A machine-oriented logic based on the resolution principle.
\newblock {\em Journal of the ACM}, 12:23--41, 1965.

\bibitem[Sho67]{qq:shoenfield}
J.~R. Shoenfield.
\newblock {\em Mathematical Logic}.
\newblock Addison-Wesley, 1967.

\end{thebibliography}
\end{document}